% Beginning polyalan.tex

\magnification = \magstep 1 \input vanilla.sty
\TagsOnRight

\title
A Generalization of P\'olya's Enumeration Theorem\\
or the  Secret Life of Certain Index Sets \\
\endtitle

\author
Valentin Vankov Iliev
\endauthor

\centerline {\it Section of Algebra, Institute of Mathematics and Informatics,}
\centerline {\it Bulgarian Academy of Sciences, 1113 Sofia, Bulgaria}
\centerline {\it E-mail: viliev\@math.bas.bg, viliev\@nws.aubg.bg}

\vskip 20pt

\heading
Introduction
\endheading

The philosophy of the present paper is that P\'olya's theory of enumeration,
which was developed in [7], can be seen as subsumed by
Schur-Macdonald's
theory of \lq\lq invariant matrices" (cf. [8], [6]).
The cornerstone of the
last one is the
equivalence between the category of finite dimensional
$K$-linear representations
of the symmetric group $S_d$ and the category of polynomial
homogeneous degree
$d$\/ functors on the category of finite-dimensional $K$-linear spaces ($K$ is
a field of characteristic $0$).
A natural background for generalization of
P\'olya's theory are the induced monomial representations of $S_d,$ which
correspond via
that equivalence to the so called semi-symmetric powers (cf. [4]). Both objects
are determined by fixing a permutation group $W\leq S_d$ and a one-dimensional
$K$-valued character  $\chi$ of $W.$ Then the representation $ind_W^{S_d}(\chi
)$ corresponds to
the $d$ th semi-symmetric power $[\chi ]^d(-).$ In particular,
when $W=S_d$ and
$\chi$ is the alternating character or the unit character, we obtain the
exterior power
$\bigwedge^d(-)$ or the symmetric power $S^d(-),$ respectively. If $W$ is the
unit group, then the regular representation of $S_d$ corresponds to the tensor
power $\otimes^d\hskip-1pt (-).$

A substantial part of problems of combinatorial analysis deal with a finite set
of objects, often called figures, and the number of figures is usually
irrelevant, provided that it is large enough. For convenience, we may suppose
that the set of figures coincides with
the set $N_0$ of non-negative integers.
The figures form
configurations, that is, elements
$(j_1,\hdots ,j_d)$ of the free monoid $Mo(N_0)$ generated by $N_0.$
Let $(x_i)_{i\in N_0}$ be a family of independent variables.
One introduces a homomorphism of monoids (weight function) $w\colon
Mo(N_0)\to Q[(x_i)_{i\in N_0}],$ where the target of $w$ is the
(commutative) algebra of polynomials in $x_0,x_1,x_2,\hdots ,$ with rational
coefficients, considered to be a monoid with respect to multiplication.
In the generic case,
one may suppose
without loss of generality
that $w$ is the
canonical homomorphism of the free monoid
onto
the free commutative monoid, both
generated by $N_0:$
$w(i)=x_i,$ $i\in N_0;$ the other cases can be
obtained by an appropriate specialization of variables $x_i.$ Under this
assumption, the monomial $x_{j_1}\hdots x_{j_d}$ is the weight of the
configuration
$(j_1,\hdots ,j_d),$ the weight function $w$ is
$S_d$-compatible and induces a bijection between the orbit space $S_d\backslash
N_0^d$ and the set of degree $d$ monomials in $x_0,x_1,x_2,\hdots .$

In general, the main
characters of P\'olya's play are the orbit spaces $W\backslash N_0^d,$
where
$W\leq S_d$
is a permutation group acting in a standard way on the integer-valued
hypercube $N_0^d\subset Mo(N_0).$
It turns out that there exist bijections of certain orbit
subspaces,
including the whole orbit space, onto
some index
sets which arise naturally within the boundaries of the semi-symmetric powers.
Given
a $K$-linear space $E$ with basis $(v_i)_{i\in N_0},$ we can construct a basis
$(v_j)_{j\in J\left (N_0^d,\chi\right )}$ for $[\chi ]^d(E).$
Here
$J(N_0^d,\chi )$ is a system of distinct representatives of
a set of $W$-orbits in $N_0^d$
and
the corresponding orbit subspace depends on the character $\chi .$ In
particular, if $\chi$
is the unit character $1_W,$ then $J(N_0^d,1_W)$ is a transversal for all
$W$-orbits in $N_0^d.$

The formal infinite sum of monomials
$$
g(W;x_0,x_1,x_2,\hdots )=
\sum_{j\in J\left (N_0^d,1_W\right )}x_{j_1}\hdots
x_{j_d}
$$
is a homogeneous degree $d$\/ symmetric function in a countable set of
variables
and it counts the $W$-orbits in $N_0^d$ provided with weights.
By the fundamental theorem on symmetric functions and by Newton's formulae
$$
g(W;x_0,x_1,x_2,\hdots )=
Z(W;p_1,\hdots ,p_d),\tag P
$$
where $p_s=\sum_{i\in N_0}x_i^s$ are the power sums and
$Z(W;p_1,\hdots ,p_d)$
is an uniquely defined isobaric polynomial in $p_1,\hdots p_d$
with rational coefficients. P\'olya's
enumeration theorem asserts that this polynomial
coincides with the cyclic index of the group $W.$ By specialization
$x_{n+1}=x_{n+2}=\cdots =0$ in (P) we establish the finite version of
P\'olya's theorem,
where the symmetric polynomial $g(W;x_0,x_1,x_2,\hdots ,x_n,0,0,\hdots )$ in
$n+1$ variables  counts the weighted $W$-orbits in the
integer-valued hypercube $[0,n]^d.$

The left-hand side of equality (P) is {\it a fortiori}
the characteristic of the polynomial functor $[1_W]^d(-)$
while
the right-hand side is
the characteristic of
the induced
monomial representation
$ind_W^{S_d}(1_W).$
The coincidence of these
characteristics is a consequence of Schur-Macdonald's equivalence.

In general, the characteristic
$ch([\chi ]^d(-))$
is a
symmetric function $g(\chi ;x_0,x_1,x_2,\hdots )$ which inventories the
weighted $W$-orbits
$O$ in $N_0^d$ subject to the condition
\lq\lq $O$ contains $|W:H|$ in number $H$-orbits," where $H\leq W$
is the kernel of $\chi$ and $|W:H|$ is the index of $H$ in $W.$
On the other hand,
the characteristic
$ch(ind_W^{S_d}(\chi ))$
is equal to the so called generalized cyclic index
$Z(\chi ;p_1,\hdots ,p_d)$
of the group $W$ with respect to the character $\chi$.
Schur-Macdonald's equivalence yields the identity
$$
g(\chi ;x_0,x_1,x_2,\hdots )=
Z(\chi;p_1,\hdots ,p_d),
$$
which turns into P\'olya's theorem in case $\chi=1_W.$

For instance, let $W$ be the whole symmetric group $S_d$ and let $\chi$ be the
alternating character $\varepsilon_d.$
Then
$g(\varepsilon ;x_0,x_1,x_2,\hdots )=ch(\bigwedge^d(-))$
is the $d$ th elementary symmetric
function $e_d$ which inventories the $S_d$-orbits
in $N_0^d$ consisting of sequences
$(i_1,\hdots ,i_d)$ with pairwise distinct components.
By a tradition dating back to P\'olya,
the formal expression $Z(A_d-S_d;p_1,\hdots ,p_d)$
is widely used
in combinatorial literature
for denoting the generalized cyclic index
$Z(\varepsilon_d;p_1,\hdots ,p_d).$
Since $\varepsilon_d=ind_{A_d}^{S_d}(1)-1_{S_d}$
in the appropriate Grothendieck group, then
$$
Z(\varepsilon_d ;p_1,\hdots ,p_d)=
Z(A_d;p_1,\hdots ,p_d)-Z(S_d;p_1,\hdots ,p_d).
$$
Thus, in the present context, the \lq\lq old" notation
$Z(A_d-S_d;p_1,\hdots ,p_d)$
can be viewed as an archetype.

In Section 1 we give a construction of a basis for the
semi-symmetric
power $[\chi ]^d(E),$ starting from a basis for the $K$-linear space $E.$
Section 2 is devoted to the main result of that paper -- Theorem 2.1.2,
which generalizes
P\'olya's Hauptsatz from
[7, Ch. 1, $n^o$ 16].
Moreover, Propositions 2.2.1 and 2.2.2  show
that the set of semi-symmetric powers is stable with respect to the tensor
product and the composition of polynomial functors (cf. [6]). The corresponding
canonical isomorphisms of functors yield generalizations of P\'olya's product
and insertion rules
[cf. 7, Ch. 1, $n^o$ 27]).
In Section 3 we give an interpretation of Theorem 2.1.2 in combinatorial terms.

\vskip 20 pt

\heading
1. Basis of a semi-symmetric power
\endheading

Throughout the rest of the paper we fix a field
$K$ of characteristic $0.$ All linear spaces and linear maps under
consideration are
assumed to be $K$-linear. Given a finite group $G,$ by a (linear)
representation of $G$ we mean a $K$-linear representation. Moreover,
if a one-dimensional character of $G$ is given, then, by default, it is
$K$-valued.
By $N_0$ we denote the set of non-negative integers. For any finite set $S$
the number of its elements is denoted by
$|S|.$

The results from this section have been announced in [3].

\vskip 17pt

1.1. Let $G$ be a finite group. We call a $G$-module any
linear representation $M$ of the group $G.$
Any element $a$ of the group ring $KG$ defines a linear endomorphism
$a\colon M\to M$ by the formula $z\mapsto az.$

Let M be a $G$-module. Given a one-dimensional character $\alpha$
of $G,$ we set:

$a_\alpha=|G|^{-1}\sum_{g\in G}\alpha (g)g,$\ $a_\alpha\in KG,$

${}_\alpha M=$\  the $G$-submodule of $M$ generated by all differences
$\alpha^{-1}(g)z-gz,$

\ \ \ \ \ \ where $g\in G$ and $z\in M,$

$M_\alpha =$\  the $G$-submodule of $M$ consisting of all $z\in M$ such
that $gz=\alpha^{-1}(g)z,$

\ \ \ \ \ \ for any $g\in G.$

Clearly, $M_\alpha$ is the isotypical component of $M,$ afforded by
the character $\alpha^{-1}.$

The following lemma is easily verified and well known
(cf. [9, Ch. I, sec. 2, $n^o$ 6, Theorem 8]).

\vskip 17pt

\proclaim {Lemma 1.1.1} Let $M$ be a $G$-module. Then

(i) The $G$-submodule ${}_\alpha M$ of $M$ is the kernel of
$a_\alpha ;$

(ii) The $G$-submodule $M_\alpha$ of $M$ is the image of
$a_\alpha .$

\endproclaim

\vskip 20pt

1.2. Let $M$ be a linear space and let $G$ be a finite group. Structure
of a {\it monomial} $G$-module on $M$ can be defined by the following data:
(a) A basis $(v_i)_{i\in I}$ for M;
(b) An action of $G$ on the index set $I;$
(c) A family $(\gamma_i)_{i\in I}$ of maps $G\to
K\backslash\{0\}$ such that
$\gamma_i(gh)=\gamma_{hi}(g)\gamma_i(h)$
for $i\in I$ and $g, h\in G.$
Then the corresponding (monomial) action of $G$
on $M$ is defined by the rule
$$
gv_i=\gamma_i(g)v_{gi}.\tag 1.2.1
$$

We write $G_i$ for the stabilizer of $i\in I$ in the group $G$ and
$G^{(i)}$ for a left transversal of $G_i$ in $G.$ Note that the
restriction of $\gamma_i$ on $G_i$ is a one-dimensional character
of the group $G_i.$

Let $\alpha$ be a one-dimensional
character of the group $G.$ Let
$I(M,\alpha )$ be the set of all $i\in I$ such that the maps
$\gamma_i$ and $\alpha^{-1}$ coincide on the subgroup $G_i.$
\vskip 17pt

\proclaim {Lemma 1.2.2} (i) The set $I(M,\alpha )$ is a $G$-stable subset of
$I.$ (ii) One has $a_\alpha (v_i)=0$ for $i\in I\backslash I(M,\alpha ).$

\endproclaim

\vskip 17pt

\demo {Proof} (i)\quad Given $i\in I,$ suppose $g\in G$ and $h\in G_i.$  Then
$G_{gi}=gG_ig^{-1}$ and $\alpha (ghg^{-1})=\alpha (h).$ Moreover,
$\gamma_{gi}(ghg^{-1})=\gamma_{g^{-1}gi}(gh)\gamma_{gi}(g^{-1})=
\gamma_{gi}(g^{-1})\gamma_i(gh)=
\gamma_{ghi}(g^{-1})\gamma_i(gh)=
\gamma_i(g^{-1}gh)=
\gamma_i(h).$

(ii)\quad The complement of $I(M,\alpha )$ in $I$ also is $G$-stable; let $i\in
I\backslash I(M,\alpha ).$
We have
$$
\align
a_\alpha (v_i) & =
|G|^{-1}\sum_{g\in G^{(i)}}\sum_{h\in G_i}\alpha (gh)\gamma_i(gh)v_{ghi} \\
& =|G|^{-1}\sum_{g\in G^{(i)}}\alpha (g)\gamma_i(g)(\sum_{h\in
G_i}\alpha (h)\gamma_i(h))v_{gi}, \\
\endalign
$$
and the equality
$a_\alpha (v_i)=0$ holds because the product $\alpha\gamma_i$ is not the unit
character of the group $G_i.$

\enddemo

\vskip 17pt

We fix a system $I^*$ of
distinct representatives of all $G$-orbits in
$I.$ Moreover, we set $J(M,\alpha )=
I^*\cap I(M,\alpha )$ and $J_0(M,\alpha
 )=I^*\backslash
J(M,\alpha ).$ Following [1, Ch. III, sec. 5, $n^o$ 4], we obtain a basis
for the linear space $M$ consisting of
$$
\align
& v_j,\quad j\in J(M,\alpha ),\tag 1.2.3 \\
& v_i-\alpha (g)\gamma_i(g)v_{gi},\quad i\in I^*,\ g\in G^{(i)},
\ g\not\in G_i, \tag 1.2.4 \\
& v_i,\quad i\in J_0(M,\alpha ).\tag 1.2.5 \\
\endalign
$$
\vskip 17pt

\proclaim {Proposition 1.2.6} Let $G$ be a finite group and let $\alpha$ be a
one-dimensional character of $G.$ Then
the following three statements hold for every
monomial $G$-module $M$ defined by the formula (1.2.1):

(i) The union of the families (1.2.4) and (1.2.5) is a basis for ${}_\alpha
M;$

(ii) The family $a_\alpha (v_j),$ \ $j\in J(M,\alpha ),$ is a basis for
$M_\alpha ;$

(iii) The family $v_j\ mod({}_\alpha M),$ \ $j\in J(M,\alpha ),$ is a
basis for the factor-space $M/{}_\alpha M.$

\endproclaim

\vskip 17pt

\demo {Proof} (i)\quad The family (1.2.4) is in ${}_\alpha M$ by definition.
Lemma
1.2.2, (ii), and Lemma 1.1.1, (i), imply that the family (1.2.5) is contained
in ${}_\alpha M.$ Now, set $J=J(M,\alpha )$ and suppose that $\sum_{j\in
J}k_ja_\alpha (v_j)=0$ for some $k_j\in K$ such that $k_j=0$ for all but
a finite number of indices $j\in J.$ We have
$$
\sum_{j\in J}k_ja_\alpha (v_j)=|G|^{-1}\sum_{j\in J}\sum_{g\in
G^{(j)}}k_j|G_j|\alpha (g)\gamma_j(g)v_{gj},
$$
hence $k_j=0$ for all $j\in J,$ which proves part (i). In
addition, we have proved that
the elements $a_\alpha (v_j),$ \ $j\in J(M,\alpha ),$ are linearly
independent.

(ii)\quad Lemma 1.1.1, (ii), yields that
the elements $a_\alpha (v_j),$ \ $j\in J(M,\alpha ),$ are in $M_\alpha $
and,
moreover, that each element of $M_\alpha$ has the form $a_\alpha (z)$ for
some $z\in
M.$ Since the union of families (1.2.3) -- (1.2.5) is a basis for $M$ and since
the endomorphism $a_\alpha$ annihilates (1.2.4) and (1.2.5), part (ii)
holds.

(iii)\quad Lemma 1.1.1 and part (ii) imply (iii).

\enddemo

\vskip 20pt

1.3. Let $W\leq S_d$ be a permutation group and let $\chi$ be a
one-dimensional character of $W.$ We recall several definitions
from [2]. Let $E$ be a linear space with basis $(v_i)_{i\in L}.$ Then
the $d$ th tensor power
$\otimes^d\hskip-1pt E$
has a standard structure of a monomial $W$-module via the rule
$$
\sigma (v_{i_1}\hskip-2pt\otimes\cdots\otimes\hskip-2pt  v_{i_d})=
v_{i_{\sigma^{-1}\left (1\right )}}\hskip-2pt \otimes\cdots\otimes\hskip-2pt
v_{i_{\sigma^{-1}\left (d\right )}}
$$
for $\sigma\in W$ and $(i_1,\hdots ,i_d)\in L^d.$

The linear factor-space $[\chi ]^d(E)=
\otimes^d\hskip-1pt E/{}_\chi (\otimes^d\hskip-1pt E)$ is called {\it $d$
\hbox{\rm th} semi-symmetric power of weight $\chi $ of $E.$} The image of the
tensor $x_1\hskip-2pt \otimes\cdots\otimes\hskip-2pt  x_d$ via the canonical
homomorphism $\otimes^d\hskip-1pt E\to [\chi ]^d(E)$ is denoted by
$x_1\chi\hdots\chi x_d.$
Clearly, the vectors of type $x_1\chi\hdots\chi x_d,$ $x_s\in E,$ generate
the semi-symmetric power $[\chi ]^d(E)$ as a linear space. The {\it $d$
\hbox{\rm th} semi-symmetric power of weight $\chi$} of a linear map $l\colon
E\to E^\prime$ is defined by the formula
$$
[\chi ]^d(l)\colon [\chi ]^d(E)\to [\chi ]^d(E^\prime ), \tag 1.3.1
$$
$$
([\chi ]^d(l))(x_1\chi\hdots\chi x_d)=
l(x_1)\chi\hdots\chi l(x_d).
$$
Varying the arguments $E$ and $l,$ we obtain
a polynomial homogeneous degree $d$\/ functor $[\chi ]^d(-)$ on the category of
(finite-dimensional) linear spaces (cf. [2], [6]).

Proposition 1.2.6 for $G=W,$ $\alpha=\chi$ and
$M=\otimes^d\hskip-1pt E$ yields the following

\vskip 17pt

\proclaim {Proposition 1.3.2}
Let $W\leq S_d$ be a permutation group and let $\chi$ be a
one-dimensional character of $W.$
Let $E$ be a linear space with basis
$(v_i)_{i\in L}.$ Then the family $(v_j)_{j\in J\left
(\otimes^d\!E,\chi\right )},$ where
$v_j=v_{j_1}\chi\hdots\chi v_{j_d},$ is a basis for the $d$ \hbox{\rm th}
semi-symmetric power $[\chi ]^d(E).$

\vskip 17pt

\proclaim {Remark 1.3.3} {\rm When the set $L$ is well-ordered, the Cartesian
product $I=L^d$ is lexicographically well-ordered. If the opposite is not
stated, we suppose that the elements of $I^*$ are lexicographically minimal
in their $W$-orbits. In this case we write $J(L^d,\chi )$ for
$J\left
(\otimes^d\hskip-1pt E,\chi\right ).$
In particular, if $L$ is the integer-valued interval $[0,n],$ where $n\in N_0,$
then we denote the index set $J(L^d,\chi )$ by $J(n,d,\chi ).$ }

\endproclaim

\vskip 20pt

\heading {2. The theorem}
\endheading

2.1. Again let $W\leq S_d$ be a
permutation group and let $\chi$ be a
one-dimensional character of $W.$
Let $x_0,x_1,x_2,\hdots$ be an infinite sequence of independent variables.
We set
$$
g_n(\chi ;x_0,\hdots ,x_n)=
\sum_{j\in J\left (n,d,\chi\right )}x_{j_1}\hdots
x_{j_d}
$$
and
$$
g(\chi ;x_0,x_1,x_2,\hdots )=
\sum_{j\in J\left (N_0^d,\chi\right )}x_{j_1}\hdots x_{j_d},
$$
where $J(n,d,\chi )$ and $J(N_0^d,\chi )$ are the index sets from (1.3.3).

\vskip 17pt

\proclaim {Lemma 2.1.1} The sequence $(g_n)_{n\geq 0}$ determines a
symmetric function
in a countable set of variables
$x_0,x_1,x_2,\hdots ,$ which coincides both with the characteristic of the
polynomial functor $[\chi ]^d(-)$ and with the formal infinite sum of
monomials
$g(\chi ;x_0,x_1,x_2,\hdots ).$

\endproclaim

\vskip 17pt

\demo {Proof} Let $v_0,\hdots ,v_n$ be a basis for the linear space
$E=K^{n+1}.$
Given an element $x=(x_0,\hdots ,x_n)\in K^{n+1}$ let $(x)$ denote the
linear
endomorphism of $E$ defined by the diagonal matrix $diag(x_0,\hdots ,x_n)$ with
respect to that basis. Then, according to (1.3.1) and (1.3.2), the linear
endomorphism $[\chi ]^d((x))$ of the semi-symmetric power $[\chi ]^d(E)$ is
defined by the diagonal matrix $diag(x_{j_1}\hdots x_{j_d})$ with respect to
the basis $(v_j)_{j\in J\left (n,d,\chi\right )}.$ In particular, the
polynomial
$g_n$ is the trace of the endomorphism $[\chi ]^d((x)).$ Therefore
$g_n$
is a symmetric polynomial in the $n+1$ variables $x_0,\hdots ,x_n$ and,
moreover, the sequence $(g_n)_{n\geq 0}$ is a
symmetric function $g(x_0,x_1,x_2,\hdots )$
in a countable set of variables,
which coincides with the characteristic of the polynomial
functor $[\chi ]^d(-)$
(cf. [6, Ch. I, Appendix, A7]).
Since
$J(n,d,\chi )=J(N_0^d,\chi )\cap [0,n]^d,$ then
$g(\chi ;x_0,\hdots x_n,0,0,\hdots )=g_n(x_0,\hdots ,x_n)$ for any $n\in N_0.$
In other words, we have
$g(\chi ;x_0,x_1,x_2,\hdots )=
g(x_0,x_1,x_2,\hdots ).$

\enddemo

\vskip 17pt

Let $p_1,\hdots ,p_d$ be independent variables. We set
$$
Z(\chi;p_1,\hdots ,p_d)=|W|^{-1}\sum_{\sigma\in W}\chi (\sigma
)p_1^{c_1\left
(\sigma\right )}\hdots p_d^{c_d\left (\sigma\right )},
$$
where $c_s(\sigma )$ is
the
number of cycles of length $s$ in the cyclic decomposition of the permutation
$\sigma\in W.$
If $\chi$ is the unit character, then
$Z(1_W;p_1,\hdots ,p_d)$
is the standard cyclic index
$Z(W;p_1,\hdots ,p_d)$
of the group $W.$

\vskip 17pt

\proclaim {Theorem 2.1.2} Let $W\leq S_d$ be a permutation group and
let $\chi$ be a
one-dimensional character of $W.$ Then one has
$$
g(\chi ;x_0,x_1,x_2,\hdots )=
Z(\chi;p_1,\hdots ,p_d), \tag 2.1.3
$$
where $p_s$ are the power sums in the variables $x_0,x_1,x_2,\hdots .$

\endproclaim

\vskip 17pt

\demo {Proof} By Frobenius reciprocity, the right-hand side of equality
(2.1.3) is the
characteristic of the induced monomial representation $ind_W^{S_d}(\chi )$ of
$S_d.$ Lemma 2.1.1 yields that the left-hand side is
the characteristic of the polynomial functor $[\chi ]^d(-).$ Due to [4],
$ind_W^{S_d}(\chi )$ and $[\chi ]^d(-)$
correspond each other via Schur-Macdonald's equivalence. Then [6, Ch. I,
Appendix, A7] implies that both
characteristics coincide.

\enddemo

\vskip 17pt

Letting
$x_{n+1}=x_{n+2}=\cdots =0$ we obtain

\vskip 17pt

\proclaim {Corollary 2.1.4} For any $n\in N_0$ one has
$$
g_n(\chi ;x_0,\hdots ,x_n)=
Z(\chi;p_1,\hdots ,p_d),
$$
where $p_s$ are the power sums in the $n+1$ variables
$x_0,\hdots ,x_n.$

\endproclaim

\vskip 17pt

\proclaim {Remark 2.1.5} {\rm If\/ $\chi$ in Theorem 2.1.2 is the unit
character,
then we establish P\'olya's classical theorem. }

\endproclaim

\proclaim {Remark 2.1.6} {\rm Corollary 2.1.4 for $n=0$ turns into the so
called \lq\lq orthogonality
relations" among the one-dimensional characters of the group $W.$
(cf. [9, Ch. I, sec. 2, $n^o$ 3, Theorem 3]).}

\endproclaim

\vskip 20pt

2.2. Now, we shall prove two statements concerning tensor product and
composition of semi-symmetric powers as polynomial functors
(cf. [6, Ch. I, Appendix, A7]).

Let $W\leq S_d$ and $V\leq S_r$ be permutation groups
and let $\chi$ and
$\theta$ be one-dimensional characters of $W$ and $V,$ respectively.
We consider the Cartesian product $S_d\times S_r$ as a subgroup of $S_{d+r}$ by
identifying $(\sigma ,\tau)\in S_d\times S_r$ with the permutation $s\mapsto
\sigma (s),$ $1\leq s\leq d,$ $d+t\mapsto d+\tau (t),$ $1\leq t\leq r.$ The
tensor product $\lambda=\chi\hskip-2pt\otimes\hskip-2pt\theta$ is a
one-dimensional character of the subgroup $W\times V\leq S_d\times S_r.$
We
identify the wreath product $S_r\hskip-2pt\sim\hskip-2pt S_d$ with its natural
faithful permutation representation in the symmetric group $S_{dr}$
(cf. [5, Ch. 4, sec. 1, 4.1.18]).
Then the wreath product $V\hskip-2pt\sim\hskip-2pt W$
is a permutation subgroup of
$S_r\hskip-2pt\sim\hskip-2pt S_d$ and the tensor product
$\mu=\theta^{\otimes d}\hskip-2pt\otimes\hskip-2pt\chi$ is
a one-dimensional character of $V\hskip-2pt\sim\hskip-2pt W.$

\vskip 17pt

\proclaim {Proposition 2.2.1} (i) For any
linear space $E$ the formulae
$$
\Pi_E\colon [\chi ]^d(E)\otimes [\theta ]^r(E)\to
[\lambda ]^{d+r}(E),
$$
$$
(y_1\chi\hdots\chi y_d)\otimes (z_1\theta\hdots\theta z_r)
\mapsto
y_1\lambda\hdots\lambda y_d\lambda z_1\lambda\hdots\lambda z_r
$$
where $y_s,\ z_t\in E,$ give rise to a canonical isomorphism of linear spaces;

(ii) The family $\Pi=(\Pi_E),$ where $E$ runs through all
finite-dimensional linear spaces, establishes a canonical isomorphism
$$
\Pi\colon [\chi ]^d(-)\otimes [\theta ]^r(-)\to
[\lambda ]^{d+r}(-)
$$
of polynomial functors;

(iii) One has
$$
g(\chi\hskip-2pt\otimes\hskip-2pt\theta ;x_0,x_1,x_2,\hdots )=
Z(\chi ;p_1,\hdots ,p_d)
Z(\theta ;p_1,\hdots ,p_r),
$$
where $p_s$ are the power sums in the variables $x_0,x_1,x_2,\hdots .$

\endproclaim

\vskip 17pt

\demo {Proof} (i)\quad This is proved in [2, sec. 2, Corollary 2.1.4].

(ii)\quad It follows directly from (1.3.1) that
the family $\Pi$ is a morphism of functors. Then part (i) implies that
$\Pi$ is an isomorphism.

(iii)\quad Since the trace is multiplicative with respect to tensor products,
then for each $n\in N_0$ we have
$
g_n(\chi\hskip-2pt\otimes\hskip-2pt\theta ;x_0,\hdots ,x_n)=
g_n(\chi ;x_0,\hdots ,x_n)
g_n(\theta ;x_0,\hdots ,x_n).
$
Therefore,
$$
g(\chi\hskip-2pt\otimes\hskip-2pt\theta ;x_0,x_1,x_2,\hdots )=
g(\chi ;x_0,x_1,x_2,\hdots )
g(\theta ;x_0,x_1,x_2,\hdots )
$$
and Theorem 2.1.2 completes the proof.

\enddemo

\vskip 17pt

\proclaim {Proposition 2.2.2} \hbox{\rm (P\'olya's insertion)} (i) For any
linear space $E$ the formulae
$$
\Delta_E\colon [\chi ]^d([\theta ]^r(E))\to
[\mu ]^{dr}(E),
$$
$$
(z_1\theta\hdots\theta z_r)\chi\hdots\chi
(z_{\left (d-1\right )r+1}\theta\hdots\theta z_{dr})\mapsto
z_1\mu\hdots\mu z_r\mu\hdots\mu
z_{\left (d-1\right )r+1}\mu\hdots\mu z_{dr},
$$
where $z_t\in E,$ give rise to a canonical isomorphism of linear spaces;

(ii) The family $\Delta=(\Delta_E),$ where $E$ runs through all
finite-dimensional linear spaces, is a canonical isomorphism
$$
\Delta\colon [\chi ]^d(-)\circ [\theta ]^r(-)\to
[\mu ]^{dr}(-)
$$
of polynomial functors;

(iii) One has
$$
g(\theta^{\otimes d}\hskip-2pt\otimes\hskip-2pt\chi ;x_0,x_1,x_2,\hdots )=
Z(\chi ;P_1,\hdots ,P_d),
$$
where
$P_s=
Z(\theta ;p_s,p_{2s},\hdots ,p_{rs}),$ $1\leq s\leq d,$
and $p_k,$ $k\geq 1,$ are the power sums in the variables $x_0,x_1,x_2,\hdots
.$

\endproclaim

\vskip 17pt

\demo {Proof} (i)\quad Due to [2, sec. 1 and sec. 2, Lemma 2.1.2], it is enough
to note that:

(a) The expression
$$
(z_1\theta\hdots\theta z_r)\chi\hdots\chi
(z_{\left (d-1\right )r+1}\theta\hdots\theta z_{dr})
\in [\chi ]^d([\theta ]^r(E))
$$
is multilinear and semi-symmetric of weight $\mu ;$

(b) The expression
$$
z_1\mu\hdots\mu z_r\mu\hdots\mu
z_{\left (d-1\right )r+1}\mu\hdots\mu z_{dr}
\in
[\mu ]^{dr}(E)
$$
is multilinear and semi-symmetric of weight $\theta$ with respect to any group
$$
Z_s=(z_{\left (d-s\right )r+1},\hdots ,z_{\left (d-s+1\right )r}),
$$
$1\leq s\leq d,$ of variables as well as semi-symmetric of weight $\chi$ with
respect to $Z_1,\hdots ,Z_d.$

(ii)\quad We take into account part (i) and (1.3.1).

(iii)\quad Part (ii),
[6, Ch. I, Appendix, (A7.3)]
and Theorem 2.1.2 yield
$$
\align
g(\theta^{\otimes d}\hskip-2pt\otimes\hskip-2pt\chi ;x_0,x_1,x_2,\hdots )& =
g(\chi ;x_0,x_1,x_2,\hdots )\circ
g(\theta ;x_0,x_1,x_2,\hdots ) \\
& =Z(\chi ;p_1,\hdots ,p_d)\circ
Z(\theta ;p_1,\hdots ,p_r),
\endalign
$$
where in the last two rows the symbol $\circ$ denotes the plethysm of symmetric
functions. Now
[6, Ch. I, sec. 8, (8.4)]
implies part (iii).

\enddemo

\vskip 20pt

\heading {3. Combinatorial interpretation}
\endheading

3.1. Let $G$ be a finite group which acts on a set $I.$
Let $\alpha$
be a one-dimensional character of the group $G$ with kernel $H\leq G.$

\vskip 17pt

\proclaim {Lemma 3.1.1} The following statements hold:

(i) If $O$ is a $G$-orbit in $I,$
then all $H$-orbits in $O$ have equal lengths;

(ii) The equalities
$\alpha_{|G_i}=1$ and $|G_i:H_i|=1$ are equivalent for any $i\in I.$

\endproclaim

\vskip 17pt

\demo {Proof}
(i)\quad Let $i\in O.$ Since $H$ is a normal subgroup of $G,$
then $\sigma H_i\sigma^{-1}\leq H$ for $\sigma\in G.$ Therefore $|H:H_{\sigma
i}|=|H:\sigma H_i\sigma^{-1}|=|H:H_i|,$ that is, each $H$-orbit in $O$ has the
same number of elements.

(ii)\quad We have $H_i=H\cap G_i$ and
$\alpha_{|G_i}=1$ is equivalent to $G_i\subset H.$

\enddemo

\vskip 17pt

Given a $G$-orbit $O$ in $I,$ we denote by
$\tau_H(O)$ the number of $H$-orbits in $O.$
If
one has $\alpha_{|G_i}=1$ for some $i\in O$ (and, hence, for all $i\in O$),
then $O$ is said to be an {\it $\alpha$-orbit.}
Then Lemma 3.1.1 and the equality
$$
|G:H||H:H_i|=|G:G_i||G_i:H_i|, \tag 3.1.2
$$
where $i\in I,$
imply the next lemma.

\vskip 17pt

\proclaim {Lemma 3.1.3} The following two statements are equivalent:

(i) The $G$-orbit $O$ is an $\alpha$-orbit;

(ii) One has $\tau_H(O)=|G:H|.$

\endproclaim

\vskip 20pt

3.2. As a consequence of Lemma 3.1.3, we establish
a bijection between the index set
$J(M,\alpha )$ from Proposition 1.2.6,
under the additional condition
$\gamma_i(g)=1$ for all $i\in I,$\  $g\in G,$
and the set of all $\alpha$-orbits.
Moreover, the equality (3.1.2) and Lemma 3.1.1, (i),
yield that $\tau_H(O)$ is a divisor of $|G:H|$ for any $G$-orbit $O$ in the
set $I.$ Hence the $G$-orbits
$O$ which satisfy the equivalent statements of Lemma 3.1.3,
contain the maximum possible number $|G:H|$ of $H$-orbits.
Thus we
obtain a combinatorial interpretation of the bases for the isomorphic linear
spaces $M/{}_\alpha M\simeq M_\alpha$ in terms of that maximum
property. In particular, Proposition 1.3.2 and the main Theorem 2.1.2 can be
restated in combinatorial terms.

\vskip 20 pt

\heading {References}
\endheading

\noindent [1] N. Bourbaki, Alg\`ebre, Chapitre III: Alg\`ebre multilineaire
(Hermann, Paris, 1948).

\noindent [2] V. V. Iliev, Semi-symmetric algebras: General constructions, J.
Algebra 148 (1992) 479 -- 496.

\noindent [3] V. V. Iliev, Semi-symmetric algebra of a free module,
C. R. Bulg. Acad. Sci. 45, $N^o$ 10 (1992) 5 -- 7.

\noindent [4] V. V. Iliev, A note on the polynomial functors corresponding to
the monomial representations of the symmetric group, J. Pure and Appl. Algebra
87 (1993) 1 -- 4.

\noindent [5] G. D. James and A. Kerber, The representation theory of the
symmetric group, {\it in} \lq\lq Encyclopedia of Mathematics and its
Applications," Vol. 16 (Addison-Wesley Publishing Company, 1981).

\noindent [6] I. G. Macdonald, Symmetric functions and Hall polynomials
(Clarendon Press, Oxford, 1979).

\noindent [7] G. P\'olya, Kombinatorische Anzahlbestimmungen f\"ur Gruppen,
Graphen und chemische Verbindungen, Acta Math. 68 (1937) 145 -- 254. English
translation:
G. P\'olya and R. C. Read, Combinatorial Enumeration of Groups,
Graphs and Chemical Compounds (Springer-Verlag New York Inc., 1987).

\noindent [8] I. Schur, \"Uber eine Klasse von Matrizen, die sich einer
gegebenen Matrix zuordnen
lassen (Dissertation, Berlin 1901), {\it in} Gesammelte Abhandlungen,
Vol. 1, 1 -- 70.

\noindent [9] J.-P. Serre, Repr\'esentations lin\'eaires des groupes finis
(Hermann, Paris, 1967).

\end